\documentclass{article}
\usepackage{amssymb,latexsym,amsmath,amsthm,mathrsfs}
\usepackage{graphicx}
\newcommand{\comment}[1]{}

\begin{document}
\title{An algebraic problem of finding four numbers given the products
of each of the numbers with the sum of the other three\footnote{Originally published as
{\em Problema algebraicum de inveniendis quatuor numeris ex datis totidem productis uniuscuiusque horum numerorum in summas trium reliquorum},
Opera Postuma \textbf{1} (1862), 282--287.
E808 in the Enestr{\"o}m index.
Translated from the Latin by Jordan Bell,
Department of Mathematics, University of Toronto, Toronto, Ontario, Canada.
Email: jordan.bell@gmail.com}}
\author{Leonhard Euler}
\date{}
\maketitle

If the four numbers which are to be found are put $v,z,y,z$, the following four equations are obtained
\begin{eqnarray*}
v(x+y+z)&=&a,\\
x(v+y+z)&=&b,\\
y(v+x+z)&=&c,\\
z(v+x+y)&=&d.
\end{eqnarray*}

By the usual rules, three unknowns can be successively eliminated from these
equations and the fourth can then be solved for. But indeed, there is no reason
why we should favor any one of these unknowns over the others
as the one to be solved for. It is appropriate that none of them be determined
by the final equation, and instead that a new unknown be introduced that
has the same relation with them all and from which the unknowns can be defined.
Thus to this end let us take the sum of the numbers to be found
\[
v+x+y+z=2t
\]
and then the above equations will turn into these
\begin{eqnarray*}
v(2t-v)=a=2tv-v^2,&\textrm{whence}&v=t-\surd(tt-a),\\
x(2t-x)=b=2tx-x^2,&&x=t-\surd(tt-b),\\
y(2t-y)=c=2ty-y^2,&&y=t-\surd(tt-c),\\
z(2t-z)=d=2tz-z^2,&&z=t-\surd(tt-d).
\end{eqnarray*}
Therefore we have now produced a solution to the extent that
from the single quantity $t$ we can readily determine all the four numbers
that are sought; hence what remains is for us to investigate this quantity
$t$, and by substituting the values in terms of $t$ just found for
$v,x,y,z$ into the equation
\[
v+x+y+z=2t,
\]
$t$ will satisfy
\[
4t-\surd(tt-a)-\surd(tt-b)-\surd(tt-c)-\surd(tt-d)=2t,
\]
from which follows
the equation
\[
2t=\surd(tt-a)+\surd(tt-b)+\surd(tt-c)+\surd(tt-d),
\]
which can indeed be made rational by the method of Newton,\footnote{Translator:
The editors of the {\em Opera omnia} refer to p. 66 of the 1707 edition of
Newton's {\em Arithmetica universalis}.}
but this would be very cumbersome. Therefore we seek to resolve this equation
in another way.

Let us put
\begin{eqnarray*}
\surd(tt-a)=p,&\textrm{and}&v=t-p,\\
\surd(tt-b)=q,&&x=t-q,\\
\surd(tt-c)=r,&&y=t-r,\\
\surd(tt-d)=s,&&z=t-s
\end{eqnarray*}
and it will be
\[
p+q+r+s=2t.
\]
Because of the irrationals $p,q,r,s$, this equation should to be
transformed into another, in which only even powers of the letters
occur. Then by doing the substitution for the letters $p,q,r,s$,
a rational equation will
be formed from which the value of the unknown $t$ may be defined.

To this end let us form the equation
\[
X^4-AX^3+BX^2-CX+D=0,
\]
whose four roots are the given quantities $a,b,c,d$.
Therefore by the nature of equations it will be
\begin{eqnarray*}
A&=&a+b+c+d,\\
B&=&ab+ac+ad+bc+bd+cd,\\
C&=&abc+abd+acd+bcd,\\
D&=&abcd.
\end{eqnarray*}
One may then put
\[
Y=tt-X \qquad \textrm{or} \qquad X=-Y+tt;
\]
by doing this substitution we shall have the equation
\[
\begin{array}{rrrrr}
Y^4&-4ttY^3&+6t^4Y^2&-4t^6Y&+t^8\\
&+AY^3&-3AttY^2&+3At^4Y&-At^6\\
&&+BY^2&-2BttY&+Bt^4\\
&&&+CY&-Ctt\\
&&&&+D
\end{array}
\Bigg\}
=0.
\]
The four roots of this equation of $Y$ will be
\[
tt-a,\quad tt-b,\quad tt-c,\quad tt-d.
\]
For this equation let us write for brevity
\[
Y^4-PY^3+QY^2-RY+S=0,
\]
so that
\begin{eqnarray*}
P&=&4tt-A,\\
Q&=&6t^4-3Att+B,\\
R&=&4t^6-3At^4+2Btt-C,\\
S&=&t^8-At^6+Bt^4-Ctt+D.
\end{eqnarray*}
Then let
\[
Y=Z^2 \qquad \textrm{or} \qquad Z=\pm \surd Y;
\]
we will have
\[
Z^8-PZ^6+QZ^4-RZ^2+S=0
\]
and the the eight roots of this equation will be the following
\[
\begin{array}{ll}
+\surd(tt-a)=+p,&-\surd(tt-a)=-p,\\
+\surd(tt-b)=+q,&-\surd(tt-b)=-q,\\
+\surd(tt-c)=+r,&-\surd(tt-c)=-r,\\
+\surd(tt-d)=+s,&-\surd(tt-d)=-s.
\end{array}
\]
One may resolve this equation of eight dimensions into two biquadratics,
the roots of the first of which are $+p,+q,+r,+s$, the second
$-p,-q,-r,-s$; let them be\footnote{Translator: If $\alpha,\beta,\gamma,\delta$
are defined by $(Z-p)(Z-q)(Z-r)(Z-s)=Z^4-\alpha Z^3+\beta Z^2-\gamma Z+\delta$,
then $(Z+p)(Z+q)(Z+r)(Z+s)=Z^4+\alpha Z^3+\beta Z^2+\gamma Z+\delta$.}
\begin{eqnarray*}
Z^4-\alpha Z^3 +\beta Z^2 -\gamma Z+\delta&=&0,\\
Z^4+\alpha Z^3+\beta Z^2 +\gamma Z+\delta&=&0, 
\end{eqnarray*}
in which by the nature of equations it will be
\begin{eqnarray*}
\alpha&=&p+q+r+s,\\
\beta&=&pq+pr+ps+qr+qs+rs,\\
\gamma&=&pqr+pqs+prs+qrs,\\
\delta&=&pqrs.
\end{eqnarray*}
Because the product of these two biquadratic equations must be equal
to the equation of eight dimensions, it will be\footnote{Translator: Multiplying
together the two biquadratics and comparing the coefficients of powers of
$Z$ in this product and in $Z^8-PZ^6+QZ^4-RZ^2+S=0$.}
\begin{eqnarray*}
P&=&\alpha^2-2\beta,\\
Q&=&\beta^2-2\alpha \gamma +2\delta,\\
R&=&\gamma^2-2\beta \delta,\\
S&=&\delta^2.
\end{eqnarray*}
And since $\alpha=p+q+r+s$, it will be
\[
\alpha=2t
\]
and hence $\alpha^2=4tt$, whence it will be
\[
\alpha^2-2\beta=4tt-2\beta=P=4tt-A,
\]
therefore
\[
\beta=\frac{A}{2}.
\]
The second equation $Q=\beta^2-2\alpha\gamma+2\delta$ will give
\[
6t^4-3Att+B=\frac{A^2}{4}-4\gamma t+2\delta
\]
or
\[
\delta=3t^4-\frac{3}{2}Att+2\gamma t-\frac{A^2}{8}+\frac{B}{2}.
\]
Indeed the third equation $R=\gamma^2-2\beta\delta$ will yield
\[
4t^6-3At^4+2Btt-C=\gamma^2-A\delta
\]
or
\[
A\delta=-4t^6+3At^4-2Btt+C+\gamma^2;
\]
this with the previous gives
\[
\begin{array}{rrrr}
4t^6&-\frac{3}{2}A^2tt&+2A\gamma t&-\frac{A^3}{8}\\
&+2Btt&&+\frac{AB}{2}\\
&&&-C
\end{array}
\bigg\}=\gamma^2.
\]
By taking the square root we will obtain\footnote{Translator: Applying
the quadratic formula to the above equation in terms of $\gamma$.}
\[
\gamma=At \pm \surd\Big(4t^6+(2B-\frac{1}{2}A^2)tt-\frac{1}{8}A^3+\frac{1}{2}AB-C\Big)\]
and thence\footnote{Translator: Using the equation
$\delta=3t^4-\frac{3}{2}Att+2\gamma t-\frac{A^2}{8}+\frac{B}{2}$.}
\[
\delta=3t^4+\frac{1}{2}Att-\frac{1}{8}A^2+\frac{1}{2}B \pm
2t \surd\Big( 4t^6+(2B-\frac{1}{2}A^2)tt-\frac{1}{8}A^3+\frac{1}{2}AB-C \Big),
\]
whose square is
\[
\begin{array}{lllll}
25t^8&+3At^6&-\frac{5}{2}A^2t^4&-\frac{5}{8}A^3tt&+\frac{1}{64}A^4\\
&&+11B&+\frac{5}{2}AB&-\frac{1}{8}A^2B\\
&&&-4C&+\frac{1}{4}B^2
\end{array}
\]
\[
\pm (12t^5+2At^3-\frac{1}{2}A^2t+2Bt)\surd\Big(4t^6+(2B-\frac{1}{2}A^2)tt
-\frac{1}{8}A^3+\frac{1}{2}AB-C\Big),
\]
should be equal to $S$, that is, to the expression $t^8-At^6+Bt^4-Ctt+D$,
from which this equation follows
\[
\begin{array}{lllll}
24t^8&+4At^6&-\frac{5}{2}A^2t^4&-\frac{5}{8}A^3tt&+\frac{1}{64}A^4\\
&&+10B&+\frac{5}{2}AB&-\frac{1}{8}A^2B\\
&&&-3C&+\frac{1}{4}B^2\\
&&&&-D
\end{array}
\]
\[
+(12t^5+2At^3-\frac{1}{2}A^2t+2Bt)\surd\Big( 4t^6+(2B-\frac{1}{2}A^2)tt
-\frac{1}{8}A^3+\frac{1}{2}AB-C \Big)=0,
\]
which rationalized gives\footnote{Translator: Move the surd term
to the right
hand side of the equation, then square both sides, then compare
coefficients of $t$. I've checked this but I don't see a clean
way to do it.}
\[
\begin{array}{lllll}
+3A^4t^8&+\frac{5}{4}A^5t^6&+\frac{3}{16}A^6t^4&+\frac{3}{256}A^7tt&+\frac{1}{4096}A^8\\
-12A^2B&-8A^3B&-\frac{27}{16}A^4B&-\frac{9}{64}A^5B&-\frac{1}{256}A^6B\\
+24AC&+7A^2C&+\frac{7}{4}A^3C&+\frac{5}{32}A^4C&+\frac{3}{128}A^4B^2\\
-48D&+12AB^2&+\frac{9}{2}A^2B^2&+\frac{9}{16}A^3B^2&-\frac{1}{32}A^4D\\
&-8AD&+5A^2D&+\frac{5}{4}A^3D&-\frac{1}{16}A^2B^3\\
&-12BC&-7ABC&-\frac{5}{4}A^2BC&+\frac{1}{4}A^2BD\\
&&-3B^3&-\frac{3}{4}AB^3&+\frac{1}{16}B^4\\
&&-20BD&-5ABD&-\frac{1}{2}B^2D\\
&&+9C^2&+\frac{5}{2}B^2C&+D^2\\
&&&+6CD
\end{array}=0.
\]
For brevity let $E=\frac{1}{4}A^2-B$ and $u=2t$; it will be
\[
\begin{array}{lllll}
+3A^2Eu^8&+2A^3Eu^6&+9A^2E^2u^4&+12AE^3uu&+4E^4\\
+6AC&+4A^2C&+28ACE&+80ADE&-32DE^2\\
-12D&+12AE^2&+80DE&+96CD&+64D^2\\
&-8AD&+36C^2&+40CE^2\\
&+12CE&+12E^3
\end{array}=0.
\]
This equation then has four positive roots and four negative roots, equal except for sign, and thus the equation can be solved as a biquadratic
equation. Moreover, the quantities $A,B,C,D$ and $E$ are known quantities
determined from the given $a,b,c,d$, since, of course,
\begin{eqnarray*}
A&=&a+b+c+d,\\
B&=&ab+ac+ad+bc+bd+cd,\\
C&=&abc+abd+acd+bcd,\\
D&=&abcd
\end{eqnarray*}
and also
\[
E=\frac{1}{4}A^2-B.
\]
For each of the values found for $u$, the sought quantities will be
\begin{eqnarray*}
v&=&\frac{u-\surd(uu-4a)}{2},\\
x&=&\frac{u-\surd(uu-4b)}{2},\\
y&=&\frac{u-\surd(uu-4c)}{2},\\
z&=&\frac{u-\surd(uu-4d)}{2}.
\end{eqnarray*}

\begin{center}
{\Large Another solution}
\end{center}

The problem can also be solved in this way. From the original equations we have
\[
\begin{array}{lll}
a-b=(v-x)(y+z),&&b-c=(x-y)(v+z),\\
a-c=(v-y)(x+z),&&b-d=(x-z)(v+y),\\
a-d=(v-z)(x+y),&&c-d=(y-z)(v+x).
\end{array}
\]
From the first and last equations we obtain
\[
v-x=\frac{a-b}{y+z},\qquad v+x=\frac{c-d}{y-z}.
\]
Let
\[
\frac{a+b-c-d}{2}=h;
\]
it will be\footnote{Translator: From the definitions of $a,b,c,d$ given
at the beginning of the paper.}
\[
h=vx-yz
\]
and by making
\[
\frac{a+b+c+d}{2}=k
\]
it will be\footnote{Translator: From the definitions of $a,b,c,d$.}
\[
k=vx+vy+vz+xy+xv+yz,
\]
and therefore
\[
k-h=2yz+(v+x)(y+z)
\]
or
\[
k-h=2yz+\frac{(c-d)(y+z)}{y-z}=c+d,
\]
therefore
\[
2yz=\frac{2dy-2cz}{y-z} \qquad \textrm{or} \qquad yyz-yzz=dy-cz.
\]
By setting $yz=t$ this equation turns into
\[
(d-t)y-(c-t)z=0.
\]
Now put
\[
dy-cz=u
\]
and it will be
\[
y=\frac{(c-t)u}{(c-d)t} \qquad \textrm{and} \qquad
z=\frac{(d-t)u}{(c-d)t}.
\]
These values substituted into $t=yz$ yield
\[
t=\frac{(c-t)(d-t)uu}{(c-d)^2tt},
\]
from which 
follows
\[
u=\frac{(c-d)t\surd t}{\surd(cd-(c+d)t+tt)},
\]
and when this is
substituted into the values found above for $y$ and $z$
we shall have
\[
\textrm{I.}\qquad y=\frac{(c-t)\surd t}{\surd(cd-(c+d)t+tt)}
\]
and
\[
\textrm{II.}\qquad z=\frac{(d-t)\surd t}{\surd(cd-(c+d)t+tt)}.
\]
Adding and subtracting,
\[
y+z=\frac{(c+d-2t)\surd t}{\surd(cd-(c+d)t+tt)}
\qquad \textrm{and} \qquad y-z=\frac{(c-d)\surd t}{\surd(cd-(c+d)t+tt)}.
\]
One then deduces that\footnote{Translator: Since $v-x=\frac{a-b}{y+z}$ and $v+x=\frac{c-d}{y-z}$.}
\[
v+x=\frac{\surd(cd-(c+d)t+tt)}{\surd t}
\qquad \textrm{and} \qquad
v-x=\frac{(a-b)\surd(cd-(c+d)t+tt)}{(c+d-2t)\surd t},
\]
whence, again adding and subtracting and putting for the sake of brevity
\[
b+c+d-a=m \qquad \textrm{and} \qquad a+c+d-b=n,
\]
the following values are produced for $v$ and $x$\footnote{Translator:
By adding and subtracting the equations for $v+x$ and $v-x$.}
\[
\textrm{III.}\qquad v=\frac{(n-2t)\surd(cd-(c+d)t+tt)}{2(c+d-2t)\surd t}
\]
and
\[
\textrm{IV.}\qquad x=\frac{(m-2t)\surd(cd-(c+d)t+tt)}{2(c+d-2t)\surd t}.
\]

But since we have found above that $h=vx-yz$, 
substituting for $v,x,y,z$ their values
and expanding yields the following equation of four dimensions for determining
the value of $t$
\[
4t(h+t)(c+d-2t)^2=(m-2t)(n-2t)(c-t)(d-t).
\]

\end{document}